\documentclass{svproc}
\usepackage{url}
\usepackage{setspace}
\usepackage{graphicx}
\usepackage[dvips]{epsfig}
\usepackage{latexsym,amsfonts,amsmath,amssymb}
\usepackage{graphicx}
\usepackage{color}
\usepackage{epstopdf}
\usepackage{makeidx}

\newtheorem{assumption}{Assumption}{\bfseries}{\itshape}
\newtheorem{condition}{Condition}{\bfseries}{\itshape}

\begin{document}

\mainmatter

\title{Turnpikes and Random Walk}
\titlerunning{Turnpikes}
\authorrunning{A.B.Piunovskiy}

\author{Alexey B. Piunovskiy}

\institute{University of Liverpool  \\ Dept. of Mathematical Sciences, Liverpool L69 7ZL, UK \\
\email{piunov@liv.ac.uk}}

\maketitle

\begin{abstract}
In this paper we revise the theory of turnpikes in discounted Markov decision processes, prove the turnpike theorem for the undiscounted model and apply the results to the specific random walk.
\end{abstract}

\section{Introduction} 

To introduce the idea of a turnpike, let us consider the discounted Markov decision process  (MDP) $\langle{\bf S},{\bf A},P,r\rangle$ with the infinite horizon, the finite state and action spaces $\bf S$ and $\bf A$, the one-step reward $r_s(a)$, and the transition probability $P_{s,j}(a)$. The optimality equation looks like
\begin{equation}\label{ape1}
W^\infty(s)=\max_{a\in {\bf A}}\left\{r_s(a)+\beta\sum_{j\in {\bf S}} P_{s,j}(a)W^\infty(j)\right\},~~s\in {\bf S},
\end{equation}
where $\beta\in(0,1)$ is the discount factor. It is well known that this equation has  a unique bounded solution which can be constructed by the successive approximations called `value iteration':
\begin{eqnarray}
W^0 && \mbox{ is an arbitrarily fixed real-valued function on } {\bf S};\nonumber\\
W^{k+1}(s)&=&U^\beta\circ W^k(s):= \max_{a\in {\bf A}}\left\{r_s(a)+\beta\sum_{j\in {\bf S}} P_{s,j}(a)W^k(j)\right\},~~s\in {\bf S},\label{ape2}\\
&& ~~~~~k=0,1,2,\ldots ; \nonumber
\end{eqnarray}
$\lim_{k\to\infty}\max_{s\in {\bf S}} |W^k(s)-W^\infty(s)|=0$. The operator $U^\beta$ is a contraction in the space of real-valued functions $W$ on $\bf S$ (which are certainly bounded) with the uniform norm $\|W\|:=\max_{s\in {\bf S}}|W(s)|$.
A control strategy is (uniformly) optimal if and only if, in each state $s\in{\bf S}$, the decision maker applies action(s) from the set
$$D^*(s):=\left\{a\in{\bf A}:~~r_s(a)+\beta\sum_{j\in {\bf S}} P_{s,j}(a)W^\infty(j)\right\}\ne\emptyset.$$
See  \cite[\S1.2.2]{apb2} or \cite[Ch.6]{apb1}. Usually, one puts $W^0(s)\equiv 0$, but it is convenient to allow $W^0$ to be arbitrary. The function $W^\infty$ certainly does not depend on $W^0$.

Since the spaces $\bf S$ and $\bf A$ are finite, the value
$$\Delta:=\min_{s\in{\bf S}}\min_{a\in{\bf A}\setminus D^*(s)}\left\{ W^\infty(s)-\left[r_s(a)+\beta\sum_{j\in {\bf S}} P_{s,j}(a)W^\infty(j)\right]\right\}
$$
is strictly positive. (Here, as usual, $\min_{a\in\emptyset}G(a):=+\infty$.)

If $\Delta=+\infty$, then all $a\in{\bf A}$ provide the maximum in (\ref{ape1}) and $D^*(s)={\bf A}$ for all $s\in{\bf S}$, and obviously, for all $k\ge 0$,
\begin{equation}\label{ape3}
\mbox{arg max}_{a\in{\bf A}}\left\{r_s(a)+\beta\sum_{j\in {\bf S}} P_{s,j}(a)W^k(j)\right\}\subset D^*(s)~\mbox{ for all } s\in{\bf S}.
\end{equation}

Suppose $\Delta<\infty$. Then, for each $s\in{\bf S}$ and each $a\in{\bf A}\setminus D^*(s)$,
$$\left[r_s(a)+\beta\sum_{j\in {\bf S}} P_{s,j}(a)W^\infty(j)\right]\le W^\infty(s)-\Delta.$$
Let us  choose $0<\varepsilon<\frac{\Delta}{2}$ and fix $K$ such that
\begin{equation}\label{ape4}
\max_{j\in{\bf S}}|W^k(j)-W^\infty(j)|<\varepsilon~~\mbox{ for all } k\ge K.
\end{equation}
Now, for each $k\ge K$, for each $s\in{\bf S}$, if $a\notin D^*(s)$ provides the maximum in (\ref{ape2}), then
\begin{eqnarray*}
W^{k+1}(s) &=& r_s(a)+\beta \sum_{j\in{\bf S}} P_{s,j}(a) W^k(j)\le  r_s(a)+\beta \sum_{j\in{\bf S}} P_{s,j}(a) W^\infty(j)+\varepsilon\\
&\le & W^\infty(s)-\Delta+\varepsilon \le W^{k+1}(s)+\varepsilon-\Delta+\varepsilon<W^{k+1}(s).
\end{eqnarray*}
The obtained contradiction shows that, for all $k\ge K$, we have inclusion (\ref{ape3}).

We have established the classical {\bf Turnpike Theorem:}

\begin{center}\begin{minipage}{10cm}
There exists $K$ such that, for all $k\ge K$, the inclusion (\ref{ape3}) holds.
\end{minipage}\end{center}

The minimal value of such $K$ is called the `turnpike integer' and  denoted as $K^*$.

This theorem appeared in \cite{apb3}; see also \cite[Thm.6.8.1]{apb1} and \cite{apb4}, where the upper bound of $K^*$ was calculated and the dependence of $K^*$ on $\beta\in(0,1)$ was discussed.

Consider the discounted MDP with the finite horizon and the fixed terminal reward $W$ which is independent of the horizon $T$:
$$W^\pi_T(\hat s):=\mathbb{E}^\pi_{\hat s}\left[\sum_{t=1}^T\beta^{t-1} r_{S(t-1)}(A(t))+\beta^T W(S(T))\right]\to\sup_\pi.$$

Here $\hat s\in{\bf S}$ is the initial state $S(0)$, $\pi$ is the control strategy (policy), $\mathbb{E}^\pi_{\hat s}$ is the mathematical expectation with respect to the strategical measure $\mathbb{P}^\pi_{\hat s}$ on $\Omega:={\bf S}\times({\bf A}\times{\bf S})^n$, and $S(t)$, $A(t)$ are the (random) controlled and action processes, i.e., the projection functions on $\Omega$. The detailed rigorous constructions can be found in \cite{apb2,apb1}.

According to the dynamic programming principle, the Bellman function 
$$W^T(s):=\sup_\pi W^\pi_T(s),~~~T=0,1,\ldots$$
satisfies equalities (\ref{ape2}), where $W^0=W$, and, for fixed $T\ge 1$, the optimal values of $A(1),A(2),\ldots,A(T)$ after observing $S(0),S(1),\ldots,S(T-1)$ are those which provide the maximum in (\ref{ape2}) at $k=T-1,T-2,\ldots,0$, $s=S(0),S(1),\ldots,S(T-1)$ correspondingly. According to the Turnpike Theorem, if $T>K^*$, then on the first steps $t=1,2,\ldots, T-K^*$,  having observed $S(t-1)$, one has to choose the actions $A(t)$ from the set $D^*(S(t-1))$. 
This observation is useful if $D^*(s)$ is a singleton for all $s\in{\bf S}$: for large time horizon $T$, the optimal feedback control on the first steps $t=1,2,\ldots,T-K^*$ is the same, independent of the terminal reward $W$ and is represented by the mapping $D^*$. This is the so called `turnpike'. In words, the Turnpike Theorem reads:
\begin{eqnarray}
&&\mbox{If the time horizon is large, on the first steps one has}\nonumber\\
&&\mbox{to control the process as if the horizon is infinite.}\label{ape9}
\end{eqnarray}
Close to the end, on the steps $T-K^*+1,T-K^*+2,\ldots,T$, the control process $A(t)$ is transient, changing as time $t$ goes on, because one has to take into account the finite remaining horizon and the terminal reward $W$. In case the set $D^*(\hat s)=\{d_1,d_2,\ldots\}$ is not a singleton for some $\hat s\in{\bf S}$, it can happen that, for each $T>K^*$, either $d_1$ or $d_2$ is not optimal at the first step in the finite horizon MDP with the initial state $\hat s$ \cite{apb3}.

It is well known that the discounting can be interpreted as the geometrical random time horizon: see \cite[\S1.2.2]{apb2}. Turnpike theorems for models with other random time horizons were studied in \cite{apb6}. In \cite{apb5}, the time horizon consists of (long enough) cycles, and during each cycle the optimal feedback strategy is of the turnpike shape.

In Section \ref{apsec2}, we  develop the turnpike theory for the undiscounted model, and in Section \ref{apsec3} we consider an example of  controlled random walk which exhibits the turnpike property.

All over the current paper, the state and action spaces ${\bf S}$ and ${\bf A}$ are finite. For the future needs, it is convenient to accept that  ${\bf S}=\{0,1,\ldots, M-1\}$, $M\ge 1$. Below, the following notations are in use:
$e:=(1,1,\ldots,1)^T\in{\mathbb R}^M$, and $I$ is the identity $M\times M$ matrix. We say that a stochastic matrix $P$ is irreducible, or aperiodic, or ergodic (i.e., irreducible and aperiodic) if the corresponding Markov chain is so. $W_s$ is the $s$-the component of the vector $W\in{\mathbb R}^M$, and $P_{s,\cdot}$ is the $s$-th row of the square $M\times M$ matrix $P$. Each function $W:{\bf S}\to{\mathbb R}$ is identified with the vector $W\in{\mathbb R}^M$, so the both notations $W(s)=W_s$ are in use. All vectors are columns; all inequalities for vectors and matrices are component-wise.  ${\bf 0}$ is the zero vector in ${\mathbb R}^M$.

\section{Description of the Undiscounted Model and the Turnpike Theorem}\label{apsec2}

In case $\beta=1$, the turnpike theory is more problematic because the operator $U^1$ is not a contraction in the uniform norm and one cannot in general expect the convergence of the sequence $\{W^k\}_{k=0}^\infty$ to a bounded function.  If, e.g., $r_s(a)>0$ for all $(s,a)\in{\bf S}\times{\bf A}$, then the sequence $\min_{s\in{\bf S}} W^k(s)$ approaches infinity as $k\to\infty$. 

Below, we consider slightly more general iterations than (\ref{ape2}).
Namely, let $\bf D$ be the finite space of `decisions', each decision leading to the (column) vector of rewards ${\cal R}(d)\in{\mathbb R}^M$ and to the stochastic $M\times M$ matrix ${\cal Q}(d)$. The triplet $\langle {\bf D},{\cal R},{\cal Q}\rangle$ will be called `the model'.

Everywhere further we assume that the introduced objects satisfy the following requirement.

\begin{assumption}\label{apas2} For every  fixed $\beta\in[0,1]$,
for each function $W:~{\bf S}\to{\mathbb R}$, there exists $\hat d\in {\bf D}$ providing the following component-wise maximum:
$${\cal R}(d)+\beta{\cal Q}(d)W\to\max_{d\in {\bf D}}.$$
\end{assumption}

The value iterations (\ref{ape2}) are generalised to 
\begin{eqnarray}
W^0 && \mbox{ is an arbitrarily fixed real-valued function on } {\bf S};\nonumber\\
W^{k+1}&=&U^\beta\circ W^k:= \max_{d\in {\bf D}}\left\{{\cal R}(d)+\beta{\cal Q}(d)W^k\right\},\label{ape5}\\
&& ~~~~~k=0,1,2,\ldots ,\nonumber
\end{eqnarray}
and the undiscounted case corresponds to $\beta=1$ in which case we omit the $\beta$ index: $U:=U^1$.

Let us describe special cases, where, like previously, $\bf A$ is a finite action space,  for each $a\in{\bf A}$, $P(a)=\left[P_{s,j}(a)\right]_{(s,j)\in{\bf S}^2}$ is a stochastic $M\times M$ matrix and $r(a)=[r_s(a)]_{s\in{\bf S}}\in{\mathbb R}^M$ is a column vector.  To put it different, let $\langle{\bf S},{\bf A},P,r\rangle$ be an MDP and accept that $\bf D$ is the set of all mappings $d:~{\bf S}\to{\bf A}$.
\begin{itemize}
\item[(i)] Standard model:
\begin{eqnarray*}
{\cal R}(d)&=& [r_s(d(s))]_{s\in{\bf S}};\\
{\cal Q}(d)&=& \left[P_{s,j}(d(s))\right]_{(s,j)\in{\bf S}^2}.
\end{eqnarray*}
In this case iterations (\ref{ape5}) represent the standard value iteration algorithm for MDP, the same as (\ref{ape2}). Assumption \ref{apas2} is obviously satisfied: the mapping $\hat d\in{\bf D}$ is built separately for all $s\in{\bf S}$. 
\item[(ii)] More general case:
\begin{eqnarray*}
{\cal R}_0(d) &=& r_0(d(0));\\
{\cal R}_{l+1}(d) &=& r_{l+1}(d(l+1))+\sum_{j=0}^l P_{l+1,j}(d(l+1)){\cal R}_j(d),\\
&&~~~~~~~~~~l=0,1,\ldots,M-2;\\
{\cal Q}_{0,j}(d) &=& P_{0,j}(d(0));\\
{\cal Q}_{l+1,j}(d) &=& \left\{\begin{array}{ll}
\displaystyle \sum_{i=0}^l\left[P_{l+1,i}(d(l+1)){\cal Q}_{i,j}(d) \right]& \mbox{ for  } j<l+1;\\
\displaystyle \sum_{i=0}^l\left[P_{l+1,i}(d(l+1)){\cal Q}_{i,j}(d)\right]+P_{l+1,j}(d(l+1)) & \mbox{ for  } j\ge l+1.\end{array}\right.\\
&&~~~~~~~~~~l=0,1,\ldots,M-2.
\end{eqnarray*}
In this case, the iterations (\ref{ape5}) are similar  to the Gauss-Seidel version of the value iteration algorithm (see \cite[\S6.3.3]{apb1}). Assumption \ref{apas2} is satisfied: the mapping $\hat d\in{\bf D}$ can be built consecutively for $s=0,1,\ldots, M-1$. This version of the general model appears in Section \ref{apsec3}.
\end{itemize}

\begin{lemma}\label{apl13}
In the version (ii) of the general model, the matrix ${\cal Q}(d)$ is stochastic for all $d\in{\bf D}$, provided the original matrix $P(a)$ is stochastic for all $a\in{\bf A}$.
\end{lemma}

\begin{definition}\label{apd1} \begin{itemize}
\item[(a)] A decision $d^{\beta *}\in{\bf D}$ is $\beta$-discounted optimal for $\beta\in(0,1)$ if 
\begin{equation}\label{ape7}
V^{\beta *}=U^\beta\circ V^{\beta *}=\max_{d\in{\bf D}}\left\{{\cal R}(d)+\beta{\cal Q}(d)V^{\beta *}\right\}={\cal R}(d^{\beta*})+\beta{\cal Q}(d^{\beta *})V^{\beta *}.
\end{equation}
Here $V^{\beta *}\in{\mathbb R}^M$ is the unique solution to the (left) equation above.
\item[(b)] A decision $d^*\in{\bf D}$ is Blackwell optimal if it is $\beta$-discounted optimal for all $\beta\in[\beta_0,1)$ for some $\beta_0\in[0,1)$.
\item[(c)] For each $d\in{\bf D}$, 
$${\cal Q}^*(d):=\lim_{n\to\infty} \frac{1}{n}\sum_{k=0}^{n-1} {\cal Q}^k(d)$$ 
is the `limiting' \cite{apb1} or `stationary' \cite{apb8} stochastic matrix; 
$$D(d):=[I-{\cal Q}(d)+{\cal Q}^*(d)]^{-1}-{\cal Q}^*(d)$$
is the deviation matrix \cite{apb8,apb1}.
\item[(d)] A decision $d^*\in{\bf D}$ is average optimal if, for all $d\in{\bf D}$,
$${\cal Q}^*(d){\cal R}(d)\le {\cal Q}^*(d^*){\cal R}(d^*).$$
${\bf D}^*$ is the set of all average optimal decisions.
\end{itemize}
\end{definition}

All the introduced objects are well defined according to  Proposition \ref{approp1}.

Let the time horizon in the model $\langle{\bf D},{\cal R},{\cal Q}\rangle$ be infinite and consider the average reward problem described as follows.
On each step $t=1,2,\ldots$, the decision $d\in{\bf D}$ leads to the transition matrix ${\cal Q}(d)$ and, if the current state is $S(t-1)=s\in{\bf S}$, then the reward equals ${\cal R}_s(d)$.   The target is to maximise the following objective
\begin{eqnarray}
&&\lim_{T\to\infty}\frac{1}{T}{\mathbb E}^d_s\left[\sum_{t=1}^{T}{\cal R}_{S(t-1)}(d)\right]=\lim_{T\to\infty}\frac{1}{T}\sum_{t=1}^{T}{\cal Q}^{t-1}_{s,\cdot}(d){\cal R}(d) =\left[{\cal Q}^*(d){\cal R}(d)\right]_s\nonumber\\
&\to&\max_{d\in{\bf D}}.\label{ape11}
\end{eqnarray}
Here $S(t)$ is the (random) state at the time moment $t=0,1,2,\ldots$ of the time-homogeneous Markov chain with the transition matrix ${\cal Q}(d)$ and the initial state $s=S(0)\in{\bf S}$; ${\mathbb E}^d_s$ is the corresponding mathematical expectation.  The average (uniformly) optimal decision is that which provides the maximum in (\ref{ape11}) simultaneously for all initial states $s\in{\bf S}$.
In case of the standard model, we have the conventional MDP \cite{apb8,apb1}, where only the stationary deterministic control strategies (policies) are considered; more general strategies (randomized, past-dependent) do not improve the gain. As will be shown (see Proposition \ref{approp1}(e)), the average optimal decision exists also in the general case.

Let the time horizon $T\in\{0,1,2,\ldots\}$  be finite and consider the total reward problem. The terminal reward $W(S(T))$ is fixed, and different decisions $d(1),d(2),\ldots,d(T)$ are allowed on the time steps $t=1,2,\ldots,T$. Now $S(\cdot)$ is the (non-homogeneous) Markov chain; the target is to maximise the total expected reward over the time horizon $T$. 
In case of the standard model, we have the conventional MDP, and it is well known that randomised decisions do not improve the objective \cite[Prop.4.4.3]{apb1}.

According to the dynamic programming principle, the Bellman function\linebreak$W^T(s)$, that is, the maximal total expected reward on the steps $1,2,\ldots,T$ including also the terminal reward $W(S(T))$, when starting from a state $S(0)=s\in{\bf S}$, satisfies equations (\ref{ape5}) with $\beta=1$, i.e.,
$$W^0=W~~\mbox{ and}~~W^{k+1}=U\circ W^k,~~~k=0,1,2,\ldots, T-1.$$
The optimal decisions $d(t)$ on the steps $t=1,2,\ldots,T$ are those which provide the maximum in (\ref{ape5}) at $k=T-1,T-2,\ldots,0$. 

In what follows, we use the `span-seminorm'
$$sp(W):=\max_{s\in{\bf S}} W(s)-\min_{s\in{\bf S}} W(s)$$
in the space of (bounded) functions $W$ on $\bf S$. The main result of the current section, Theorem \ref{apt2}, states that the condition $\lim_{k\to\infty} sp(W^{k+1}-W^k)=0$ is sufficient for the turnpike property of the model. Note that in general the sequence $sp(W^{k+1}-W^k)$, $k=0,1,2,\ldots$ does not approach zero \cite[\S4.2.18]{apb7}, \cite[Ex.8.5.1]{apb1}; see also Example in Subsection \ref{apsec32}.

\begin{theorem}\label{apt2} {\bf (Turnpike)}
Suppose for  $W^0=W$ the iterations (\ref{ape5}) with $\beta=1$ satisfy the requirement: $\lim_{k\to\infty} sp(W^{k+1}-W^k)=0$. Then there exists $K$ such that, in the described model with the time horizon  $T> K$, on the first steps $t=1,2,\ldots, T-K$ the optimal decisions are necessarily average optimal ones (in the sense of Definition \ref{apd1}(d)).
\end{theorem}

Sufficient conditions for the requirement $\lim_{k\to\infty} sp(W^{k+1}-W^k)=0$ are provided in Section \ref{apsec6}. If this requirement is not satisfied, the turnpike property may be not valid: see Example in Subsection \ref{apsec32}.

Like in the discounted model, the Turnpike Theorem is useful if ${\bf D}^*$, the set of average optimal decisions, is a singleton. Otherwise, if, e.g., $d_1,d_2\in{\bf D}^*$ then it can happen that, for each $T$, either $d_1$ or $d_2$ is not optimal at the first step in the finite horizon model.  The situation when ${\bf D}^*$ is a singleton is described in the following statement.

\begin{theorem}\label{apco1}
Suppose for each $d\in{\bf D}$ the matrix ${\cal Q}(d)$ is irreducible, the Blackwell optimal strategy $d_B$ is unique, and $d^*$ is the unique decision providing the maximum in the optimality equation
\begin{equation}\label{ape10}
Ge+Y=\max_{d\in{\bf D}}\left\{{\cal R}(d)+{\cal Q}(d)Y\right\}.
\end{equation}
(Under the imposed conditions, the vector $Y$ is unique up to the additive constant.)
Then $d^*=d_B$ is the unique average optimal decision.
\end{theorem}

\section{Proposition \ref{approp1} and Proof of Theorems \ref{apt2} and \ref{apco1}}\label{apsec5}

\begin{proposition}\label{approp1}
\begin{itemize}
\item[(a)] A $\beta$-discounted optimal  decision exists for each $\beta\in(0,1)$.
\item[(b)] Suppose $\beta\in(0,1)$ is fixed. For each $d\in{\bf D}$, the unique solution to the equation
\begin{equation}\label{ape8}
V^{\beta d}={\cal R}(d)+\beta{\cal Q}(d)V^{\beta d}
\end{equation}
satisfies the inequality $V^{\beta d}\le V^{\beta *}$.
\item[(c)] For each $d\in{\bf D}$ the limiting matrix ${\cal Q}^*(d)$ exists, is stochastic, and exhibits the following properties:
$${\cal Q}^*(d)={\cal Q}(d){\cal Q}^*(d)={\cal Q}^*(d){\cal Q}(d)={\cal Q}^*(d){\cal Q}^*(d).$$
\item[(d)] The deviation matrix $D(d)$ is well defined and satisfies equation
$$D(d)=[I-{\cal Q}(d)+{\cal Q}^*(d)]^{-1}](I-{\cal Q}^*(d)).$$
\item[(e)] There exists a Blackwell optimal decision, and 
every Blackwell optimal decision is also average optimal. Hence an average optimal decision also exists and ${\bf D}^*\ne\emptyset$.
\end{itemize}
\end{proposition}

Proofs of all auxiliary statements are presented in Appendix.

\begin{lemma}\label{apl1}
Suppose a decision $\hat d\in{\bf D}$ provides maximum to 
$${\cal R}(d)+{\cal Q}(d)W,$$
where $W\in{\mathbb R}^M$. Then
\begin{eqnarray*}
&&\min_{s\in{\bf S}}\left\{U\circ W(s)-W(s)\right\}e\le{\cal Q}^*(\hat d) {\cal R}(\hat d)\le {\cal Q}^*(d^*) {\cal R}(d^*)\\
&\le& \max_{s\in{\bf S}}\left\{U\circ W(s)-W(s)\right\}e,
\end{eqnarray*}
where $d^*$ is an average optimal decision.
\end{lemma}

\underline{Proof of Theorem \ref{apt2}.}  Denote ${\bf D}^*$ the set of all average optimal decisions in the model $\langle {\bf D},{\cal R},{\cal Q}\rangle$. This set is not empty due to Proposition \ref{approp1}(e). Let $g^*:={\cal Q}^*(d^*){\cal R}(d^*)$ for $d^*\in{\bf D}^*$. Clearly, the vector $g^*$ does not depend on the choice of $d^*\in{\bf D}^*$. Let
$$\varepsilon:=\min_{d\in{\bf D}\setminus{\bf D}^*}\|g^*-{\cal Q}^*(d){\cal R}(d)\|>0,$$
where, as usual, $\|V\|:=\max_{s\in{\bf S}}|V(s)|$ is the uniform norm in  ${\mathbb R}^M$. Take $K$ such that
$$sp(W^{k+1}-W^k)=sp(U\circ W^k-W^k)<\varepsilon$$
for all $k\ge K$.

Now, for each $k\ge K$, if $\hat d\in{\bf D}$ provides the maximum to ${\cal R}(d)+{\cal Q}(d) W^k$ then, by Lemma \ref{apl1},
\begin{eqnarray*}
&&g^*-{\cal Q}^*(\hat d){\cal R}(\hat d)={\cal Q}^*(d^*){\cal R}(d^*)-{\cal Q}^*(\hat d){\cal R}(\hat d)\\
&\le &  sp(U\circ W^k-W^k)e<\varepsilon e
\Longrightarrow  \|g^*-{\cal Q}^*(\hat d){\cal R}(\hat d)\|<\varepsilon;
\end{eqnarray*}
hence $\hat d\in{\bf D}^*$.

Therefore, for each $T>K$, on the first steps $t=1,2,\ldots, T-K$, the optimal decisions, which provide the maximum to ${\cal R}(d)+{\cal Q}(d) W^k$ at $k=T-1,T-2,\ldots,K$, belong to ${\bf D}^*$.\hfill$\Box$\vspace{3mm}

\underline{Proof of Theorem  \ref{apco1}.}
If the matrix ${\cal Q}(d)$ is irreducible, then the stochastic matrix ${\cal Q}^*(d)$ has identical strictly positive rows: see Thm.5.1.1 and Thm.5.1.4 of  \cite{apb9}. Therefore, ${\cal Q}^*(d){\cal R}( d)=G(d)e$ for some $G(d)\in{\mathbb R}$ and $d^*\in{\bf D}^*$ if and only if $G(d^*)=\max_{d\in{\bf D}} G(d)=:G^*$. Recall also that a Blackwell optimal decision exists according to Proposition \ref{approp1}(e). 

The optimality equation (\ref{ape10})
is solvable, and a pair $(G,Y)\in{\mathbb R}\times{\mathbb R}^M$ satisfies it if and only if 
$$G=G^*~\mbox{ and } Y=ce+D(d_B){\cal R}(d_B),$$
where $c\in{\mathbb R}$ is an arbitrary constant, $d_B$ is some Blackwell optimal decision and $D(d_B)$
is the  deviation matrix.  In the case of the standard model, this is precisely Theorem 6.1 of \cite{apb8}; all the steps of its proof remain correct for the general case. 

Let us show that a decision $d^*$ is average optimal if and only if it provides the maximum in (\ref{ape10}).
Let $(G^*,Y)$ be a solution to the equation (\ref{ape10}). Then
$${\cal G}(d):={\cal R}(d)-G^*e+{\cal Q}(d)Y-Y\le {\bf 0}~\mbox{for all } d\in{\bf D},~\mbox{ and } \max_{d\in{\bf D}} {\cal G}(d)={\bf 0}.$$

Since all the rows of each matrix ${\cal Q}^*(d)$ are strictly positive,
$${\cal Q}^*(d){\cal G}(d) \le {\bf 0}~\mbox{for all } d\in{\bf D}~\mbox{ and } \max_{d\in{\bf D}} {\cal Q}^*(d) {\cal G}(d)={\bf 0}.$$
Moreover, for a decision $\hat d$ providing the last maximum, we have
\begin{equation}\label{ape24}
{\cal Q}^*(\hat d) {\cal G}(\hat d)=\max_{d\in{\bf D}} {\cal Q}^*(d) {\cal G}(d)={\bf 0}\Longleftrightarrow {\cal G}(\hat d)=\max_{d\in{\bf D}} {\cal G}(d)={\bf 0}
\end{equation}
because ${\cal Q}^*(d)>0$ for all $d\in{\bf D}$.
But ${\cal Q}^*(d) {\cal G}(d)={\cal Q}^*(d){\cal R}(d)-G^*e$ for all $d\in{\bf D}$ (see Proposition \ref{approp1}(c)), so that a decision $d^*$ is average optimal if and only if
$${\bf 0}={\cal Q}^*(d^*) {\cal G}(d^*)=\max_{d\in{\bf D}} {\cal Q}^*(d) {\cal G}(d)\Longleftrightarrow {\cal G}(d^*)=\max_{d\in{\bf D}} {\cal G}(d)\mbox{ by (\ref{ape24})},$$
i.e., if and only if $d^*$ provides the maximum in (\ref{ape10}).

Now the statement of Theorem \ref{apco1} follows. Recall that the decision $d_B$ is average optimal by Proposition \ref{approp1}(e), so that $d^*=d_B$.\hfill$\Box$

\section{Sufficient Conditions for the Equation $\lim_{k\to\infty}sp(W^{k+1}-W^k)=0$ in Case $\beta=1$}\label{apsec6}

For fixed $J\ge 1$ take arbitrary sequences ${d^1}=\{d^1_1,d^1_2,\ldots,d^1_J\}$ and ${d^2}=\{d^2_1,d^2_2,\ldots,d^2_J\}$ of decisions from $\bf D$ and  matrices
\begin{eqnarray*}
{\cal Q}^{(J,1)}&:=&{\cal Q}(d^1_1){\cal Q}(d^1_2)\ldots {\cal Q}(d^1_J);\\
{\cal Q}^{(J,2)}&:=&{\cal Q}(d^2_1){\cal Q}(d^2_2)\ldots {\cal Q}(d^2_J).
\end{eqnarray*}
Denote
$$\eta({d^1},{d^2}):= \min_{(s,u)\in{\bf S}^2} \sum_{l\in {\bf S}} \min \left\{{\cal Q}^{(J,1)}_{s,l}, {\cal Q}^{(J,2)}_{u,l}\right\}$$
and introduce
\begin{equation}\label{ape6}
\gamma^J:=1-\min_{\left({d^1},{d^2}\right)} \eta({d^1},{d^2})\in[0,1].
\end{equation}

\begin{lemma}\label{apt1}
Suppose 
$\gamma^J<1$ for some $J\ge 1$. 
Then the following statements hold.
\begin{itemize}
\item[(a)] The operator $U:=U^1$ in (\ref{ape5}) is a $J$-step contraction with respect to the span-seminorm:
$$sp(U^J\circ V_1-U^J\circ V_2)\le \gamma^J sp(V_1-V_2)$$
for all $V_1,V_2\in{\mathbb R}^M$.
\item[(b)] For every $W^0\in{\mathbb R}^M$, given $\varepsilon>0$, there exists $K$ such that for all $k\ge K$
$$sp(W^{k+1}-W^k)\le\varepsilon,$$
that is, $\lim_{k\to\infty} sp(W^{k+1}-W^k)=0$.
\end{itemize}
\end{lemma}

Sufficient conditions for the inequality $\gamma^J<1$ at some $J\ge 1$ are given in the following lemma.

\begin{lemma}\label{apt3}
Suppose there exist $J\ge 1$ 
and  a state $l\in{\bf S}$ such that, for any sequence $\{d_1,d_2,\ldots,d_J\}$ of decisions from $\bf D$, for the matrix 
$${\cal Q}^{(J)}:={\cal Q}(d_1){\cal Q}(d_2)\ldots{\cal Q}(d_J),$$
the strict inequality ${\cal Q}^{(J)}_{s,l}>0$ is valid  for all $s\in{\bf S}$.

Then $\gamma^J<1$, and hence $\lim_{k\to\infty} sp(W^{k+1}-W^k)=0$ for every $W^0\in{\mathbb R}^M$.
\end{lemma}

The proof is identical to the proof of Theorem 8.5.3(b) of \cite{apb1}.

\begin{corollary}\label{apco2}
If, for each $d\in{\bf D}$, the matrix ${\cal Q}(d)$ is ergodic  and ${\cal Q}_{s,s}(d)>0$ for all $s\in{\bf S}$, then $\gamma^J<1$ for some $J\ge 1$, and hence $\lim_{k\to\infty} sp(W^{k+1}-W^k)=0$ for every $W^0\in{\mathbb R}^M$.
\end{corollary}

\begin{lemma}\label{apt4}
Suppose for every average optimal decision $d^*\in{\bf D}^*$ the matrix ${\cal Q}(d^*)$ is aperiodic. Assume also that ${\cal Q}^*(d^*){\cal R}(d^*)=G^*e$: the maximal value of the objective (\ref{ape11}) does no depend on the initial state $s$. Then, for every $W^0\in{\mathbb R}^M$, $\lim_{k\to\infty} sp(W^{k+1}-W^k)=0$.
\end{lemma}

\section{Controlled Random Walk}\label{apsec3}

\subsection{Description of the Model and the Associated Conjecture}\label{apsec31}

The problem under study is a specific Markov Decision Process (MDP) defined as follows.

Suppose
$${\bf A}=\{a_1,a_2,\ldots, a_N\}$$
is the finite action space and, for each $a\in{\bf A}$, let $Z(a)$ be the random variable taking values $m=1,2,\ldots, M$ with probabilities $p_m(a)$. The state space of MDP is
$${\bf X}:=\{-M,-M+1,-M+2,\ldots\}.$$
All the states $-M,-M+1,\ldots, -1$ are absorbing, with zero rewards.

If action $a\in{\bf A}$ is chosen in the state $i\ge 0$, then the new state is just $j=i-Z(a)$, that is, the transition probability is given by
$$P_{i,j}(a)=\left\{\begin{array}{ll}
p_m(a), & \mbox{ if } j=i-m,~~m=1,2,\ldots, M; \\
0 & \mbox{ otherwise.}
\end{array}\right.$$
Finally, the associated (expected) reward in states $i\ge 0$ equals $R^a$ and is $i$-independent. For example, if $R_{Z(a)}(a)$ is the reward associated with the action $a\in{\bf A}$ and the value $Z(a)$, then
$$R^a=\sum_{m=1}^M R_m(a) p_m(a).$$
To summarise, the one-step rewards in the MDP equal
$$r_i(a):=\left\{\begin{array}{ll}
0, & \mbox{ if } i < 0; \\ R^a, & \mbox{ if } i\ge 0. \end{array}\right. $$
The initial state $i\in{\bf X}$ is fixed, and we consider MDP $\langle{\bf X},{\bf A},P,r\rangle$ with the total expected reward, with (random) states and actions
$$X_0=i,~A_1,~X_1,~A_2,\ldots~.$$

The definition of a strategy $\pi$ (past-dependent, randomized) is conventional \cite{apb8,apb2,apb1}; $E^\pi_i$ is the corresponding mathematical expectation;
$$V(i):=\sup_\pi E^\pi_i\left[\sum_{t=1}^\infty r_{X_{t-1}}(A_t)\right]$$
is the Bellman function for this MDP; $i\in{\bf X}$. Since the reward $r$ is bounded and the process $X_t$ is ultimately absorbed after (maximum) $i+1$ time steps  at $\{-M,-M+1,\ldots,-1\}$, the function $V$ is finite-valued. It is well known (see, e.g., \cite[Ch.4]{apb8} or \cite[\S9.5]{apb9})  that the function $V$ is the unique solution to the optimality (Bellman) equation
\begin{eqnarray}\label{ape21}
V(i) &=& \max_{a\in{\bf A}}\left\{ R^a+\sum_{m=1}^M V(i-m)p_m(a)\right\}~~~\mbox{ for } i\ge 0;\\
V(i) &=& 0 ~~~\mbox{ for } i=-M,-M+1,\ldots, -1, \nonumber
\end{eqnarray}
which can be solved successively for $i=0,1,\ldots$.

Let us introduce the following notations: $\displaystyle L^a:= \sum_{m=1}^M m p_m(a)\ge 1;~~ c_*:= \max_{a\in{\bf A}} \frac{R^a}{L^a}$, and ${\bf A}_*$ is the set of all actions providing this maximum.

In the current section, we show that, under appropriate conditions, the following {\bf conjecture}:
\begin{eqnarray}\label{ape22}
&\mbox{There is such $I<\infty$ that, for all $i \ge I$,}\\
&\mbox{ the maximum in (\ref{ape21}) is only provided by $a\in {\bf A}_*$}\nonumber
\end{eqnarray}
is valid: see Theorem \ref{apt6}. 
This conjecture was formulated by Prof.I.Sonin in a private conversation.
Example in Subsection \ref{apsec32} shows that  this conjecture may be not valid in some situations.

\subsection{Proof of the Conjecture}\label{apsec33}

We will reformulate the conjecture (\ref{ape22}) as the Turnpike Theorem provided in Section \ref{apsec2}. The first step is described in the following lemma.

\begin{lemma}\label{apl8}
The function
\begin{equation}\label{ape25}
\tilde W(i):= V(i)-c_*i,~~~i\in{\bf X}
\end{equation}
is the (unique) uniformly bounded function satisfying equation
\begin{eqnarray}
\label{ape23}
\tilde W(i)&=&-c_*i~~~ \mbox{ for } i=-M,-M+1,\ldots, -1; \nonumber\\
\tilde W(i)&=&\max_{a\in{\bf A}}\left\{ L^a\left(\frac{R^a}{L^a}-c_*\right) +\sum_{m=1}^M \tilde W(i-m) p_m(a)\right\}~~~\mbox{ for } i\ge 0,
\end{eqnarray}
which can be solved successively for $i=0,1,\ldots$. Hence $V(i)=c_*i+O(1)$ when $i\to\infty$.

Moreover, for each $i\in{\bf X}$, the maxima in (\ref{ape21}) and in (\ref{ape23}) are provided by the same values of $a\in{\bf A}$.
\end{lemma}

Every value of $i\in{\bf X}$ can be uniquely represented as
$$i=(k-1)M+s,~~\mbox{ where } s\in{\bf S}:=\{0,1,\ldots, M-1\},~~ k=0,1,2,\ldots .$$
For each $i\in{\bf X}$ with the corresponding values of $k$ and $s$, we denote  $\tilde W(i)$, introduced in (\ref{ape25}), as $W^k(s)$. 
Now equation (\ref{ape23}) takes the following form:
$$W^0(s)=-c_*(-M+s)~~~\mbox{ for } s\in{\bf S};$$ {\small
\begin{eqnarray*}
W^{k+1}(0)&=&\max_{a\in{\bf A}}\left\{ L^a\left(\frac{R^a}{L^a}-c_*\right)+\sum_{j=0}^{M-1} W^k(j) p_{M-j}(a)\right\},\\
W^{k+1}(1)&=&\max_{a\in{\bf A}}\left\{ L^a\left(\frac{R^a}{L^a}-c_*\right)+W^{k+1}(0)p_1(a)+\sum_{j=1}^{M-1} W^k(j) p_{M-j+1}(a)\right\},\\
&\ldots &\\
W^{k+1}(M-1)&=&\max_{a\in{\bf A}}\left\{ L^a\left(\frac{R^a}{L^a}-c_*\right)+W^{k+1}(M-2)p_1(a)+W^{k+1}(M-3) p_2(a)\right.\\
&&+\ldots\left.\vphantom{\frac{R^a}{L^a}}+W^{k+1}(0) p_{M-1}(a)+W^k(M-1)p_M(a)\right\}, ~~k=0,1,\ldots.
\end{eqnarray*}}
After we introduce the stochastic matrix  {\small
\begin{equation}\label{ape26}
  P(a):= \left(
\begin{array}{llll}
  P_{0,0}(a)=p_M(a) &   P_{0,1}(a)=p_{M-1}(a) & \ldots &   P_{0,M-1}(a)=p_1(a) \\
  P_{1,0}(a)=p_1(a) &   P_{1,1}(a)=p_{M}(a) & \ldots &   P_{1,M-1}(a)=p_2(a) \\
\ldots && \ldots \\
  P_{M-1,0}(a)=p_{M-1}(a) &   P_{M-1,1}(a)=p_{M-2}(a) & \ldots &   P_{M-1,M-1}(a)=p_M(a) \\
\end{array} \right),
\end{equation}}
the obtained equations for $W^k(s)$ can be rewritten as {\small
\begin{eqnarray}
W^0(s) &=& -c_*(-M+s)~~~\mbox{ for } s\in{\bf S};\label{ape29}\\
W^{k+1}(s) &=& \max_{a\in{\bf A}}\left\{ L^a\left(\frac{R^a}{L^a}-c_*\right)+\sum_{j=0}^{s-1} W^{k+1}(j)   P_{s,j}(a)+\sum_{j=s}^{M-1} W^k(j)   P_{s,j}(a)\right\}\nonumber\\
&&~~\mbox{ for } s\in{\bf S},~k\ge 0.\nonumber
\end{eqnarray}}

Let us put $r_s(a):=L^a\left(\frac{R^a}{L^a}-c_*\right)\le 0$ and consecutively for $s=0,1,2,\ldots,$\linebreak$M-1$ express $W^{k+1}(s)$ in terms of $W^k$. As the result, we finish with the equations
\begin{eqnarray}
W^0(s) &=& -c_*(-M+s),~~~s\in{\bf S};\nonumber\\
W^{k+1} &=& \max_{d\in{\bf D}}\{{\cal R}(d)+{\cal Q}(d)W^k\},~~k=0,1,\ldots, \label{ape27}
\end{eqnarray}
where $\bf D$ is the set of all mappings $d:~{\bf S}\to{\bf A}$, and ${\cal R}(d)\le {\bf 0}$ and ${\cal Q}(d)$ are as in the version (ii) of the general model $\langle{\bf D},{\cal R},{\cal Q}\rangle$ described in Section \ref{apsec2}. 

\begin{definition}\label{apd3}
Decisions $d\in{\bf D}$ satisfying the property $d(s)\in{\bf A}_*$ for all $s\in{\bf S}$ will be called trivial. Equivalently, a decision $d\in{\bf D}$ is trivial if and only if ${\cal R}(d)={\bf 0}$.
\end{definition}

The conjecture (\ref{ape22}) is now reformulated as follows:
\begin{eqnarray}\label{ape28}
&\mbox{There exists $K$ such  that, for all $k \ge K$,
 the maximum in (\ref{ape27}) }\\
&\mbox{is only provided by the trivial decisions.}\nonumber
\end{eqnarray}
Note that all the vectors $W^0,W^1,\ldots$ are uniformly bounded by Lemma \ref{apl8}, and the maxima in (\ref{ape21}),(\ref{ape23}) and  (\ref{ape29}) are provided by the same values of $a\in{\bf A}$.

Since ${\cal R}(d)\le {\bf 0}$ for all $d\in{\bf D}$ and ${\cal R}(d)={\bf 0}$ for each trivial decision $d$, it is obvious that ${\bf D}^*$, the set of average optimal decisions in the model $\langle{\bf D},{\cal R},{\cal Q}\rangle$, contains all trivial decisions: the optimal gain $\max_{d\in{\bf D}}{\cal Q}^*(d){\cal R}(d)={\bf 0}$  is attained at any trivial decision.

Suppose the iterations (\ref{ape27}) satisfy the requirement $\lim_{k\to\infty} sp(W^{k+1}-W^k)=0$. Then, according to the Turnpike Theorem \ref{apt2}, there exists $K$ such that for each $T>K$ the maximum in (\ref{ape27}) at $k=T-1$ is necessarily provided by an average optimal decision $d\in{\bf D}^*$ in the model $\langle{\bf D},{\cal R},{\cal Q}\rangle$. If ${\bf D}^*$ contains only the trivial decisions, then the conjecture (\ref{ape28}) (and also (\ref{ape22})) is valid. This observation makes possible to formulate the sufficient conditions for the conjecture (\ref{ape22}) in terms of the matrices ${\cal Q}(d)$ introduced above and based on the matrix (\ref{ape26}).

\begin{condition}\label{apcon2} For each $d^*\in{\bf D}^*$ the matrix ${\cal Q}(d^*)$ corresponds to the Markov chain without transient states.
\end{condition}

\begin{lemma}\label{apl11}
If Condition \ref{apcon2} is satisfied then  ${\bf D}^*$ contains only trivial decisions.
\end{lemma}

\begin{condition}\label{apcon3} Either 
\begin{itemize}
\item[$\bullet$] $\gamma^J<1$ for some $J\ge 1$, where $\gamma^J$ is as defined at the beginning of Section \ref{apsec6},
\item[$\bullet$] or, for each $d^*\in{\bf D}^*$, the matrix ${\cal Q}(d^*)$ is aperiodic.
\end{itemize}
\end{condition}

\begin{theorem}\label{apt6} If Conditions \ref{apcon2} and \ref{apcon3} are satisfied, then Conjecture (\ref{ape22}) is valid.
\end{theorem}

The proof directly follows from Lemmas \ref{apt1}, \ref{apt4}, Theorem \ref{apt2} and Lemma \ref{apl11}. See also the explanations above.

In the following corollary, we provide the sufficient conditions for Conjecture (\ref{ape22}) to be valid, in terms of the original random walk.

\begin{corollary}\label{apcor3} Suppose $p_M(a)>0$ for all $a\in{\bf A}$ and, for any two states $i,j\in{\bf S}$, there exists a path $i_0=i\to i_1\to\ldots\to i_N=j$ in $\bf S$ such that, for any $a_0,a_1,\ldots,a_{N-1}\in{\bf A}$, 
$$P_{i_0,i_1}(a_0)P_{i_1,i_2}(a_1)\ldots P_{i_{N-1},i_N}(a_{N-1})>0,$$
where $P(a)$ is given by (\ref{ape26}).

Then Conjecture (\ref{ape22}) is valid.
\end{corollary}

The matrix $P(a)$ has a cyclic structure. Thus, the conditions of Corollary \ref{apcor3} are satisfied if there is $m<M$ having no common divisors with $M$ such that $p_M(a)>0$ and $p_m(a)>0$ for all $a\in{\bf A}$.

\subsection{Example}\label{apsec32}

In this subsection, we show that the conjecture (\ref{ape22}) may be not valid if  Conditions \ref{apcon2} and \ref{apcon3}  are not satisfied. We use the notations from Subsection \ref{apsec31}. In particular, the function $V$ comes from the iterations (\ref{ape21}).

Put 
$${\bf A}:=\{a_1,a_2\},~M:=3,~\varepsilon\in(0,1),~  p_2(a_1)=1,~ p_2(a_2)=1-\varepsilon,~p_3(a_2)=\varepsilon,$$
where $\varepsilon\in(0,1)$;
other probabilities being zero. Finally, let $R^{a_1}:=2$ and $R^{a_2}:=h\in(2,2+\varepsilon)$. Now
$$L^{a_1}=2,~L^{a_2}=2+\varepsilon,~ \frac{R^{a_1}}{L^{a_1}}=1,~ \frac{R^{a_2}}{L^{a_2}}=\frac{h}{2+\varepsilon}<1,~c_*=1,~{\bf A}_*=\{a_1\}.$$

Since $h>2$, obvious calculations lead to the following expressions:\\
$V(0)=V(1)=\max\{2,h\}=h;$\\
$V(2)=\max\{2+V(0)=2+h;~~h+(1-\varepsilon)V(0)+\varepsilon V(-1)=h+(1-\varepsilon)h\}=2+h$\\
because 
$$\frac{2}{1-\varepsilon}=2[1+\varepsilon+\varepsilon^2+\ldots]>2+\varepsilon>h\Longrightarrow 2>(1-\varepsilon)h.$$
$V(3)=\max\{2+V(1)=2+h;~~h+(1-\varepsilon)V(1)+\varepsilon V(0)=2h\}=2h.$
Further properties of the function $V$ are given in the following lemma.

\begin{lemma}\label{apl21}
For all $j\ge 1$, the following statements hold.
\begin{itemize}
\item[(i)] For even steps $i=2j$, 
$$V(2j)=2j+h,$$ 
and maximum in (\ref{ape21}) is provided by $a_1$ only.
\item[(ii)] For odd steps $i=2j-1$, 
$$V(2j-1)<\frac{2\varepsilon(j-1)+(1+\varepsilon)h-2}{\varepsilon}.$$
\item[(iii)] For odd steps $i=2j+1$, 
$$V(2j+1)=(1-\varepsilon)V(2j-1)+(1+\varepsilon)h+2\varepsilon(j-1),$$ 
and maximum in (\ref{ape21}) is provided by $a_2$ only.
\end{itemize}
\end{lemma}

Therefore, for all odd values of $i$, the maximum in (\ref{ape21}) is provided only by $a_2\notin {\bf A}_*$. The conjecture (\ref{ape22}) is not valid.

Now we look at this example from the viewpoint of the generalised model $\langle{\bf D},{\cal R},{\cal Q}\rangle$ described in Subsection \ref{apsec33}: $M=3$ and ${\bf S}=\{0,1,2\}$. The only trivial decision is $d^*(s)\equiv a_1$, which is certainly average optimal, as mentioned below Definition \ref{apd3}.

\begin{lemma}\label{apl22}
The only average optimal decision in the model $\langle{\bf D},{\cal R},{\cal Q}\rangle$ is the trivial decision $d^*$.
\end{lemma}

The functions $W^k$, $k\ge 1$, coming from $V$, have the following form (see Lemma \ref{apl21}):
\begin{itemize}
\item[$\bullet$] If $k$ is odd, then
\begin{eqnarray*}
W^k(0) &=& \tilde W(3(k-1))=V(3(k-1))-3(k-1)=h;\\
W^k(1) &=& \tilde W((3(k-1)+1)=V(3(k-1)+1)-[3(k-1)+1]\\
&=:& z(2j+1),~\mbox{ where } j:=3(k-1)/2;\\
W^k(2) &=& \tilde W((3(k-1)+2)=V(3(k-1)+2)-[3(k-1)+2]=h.
\end{eqnarray*}
\item[$\bullet$] If $k$ is even, then
\begin{eqnarray*}
W^k(0) &=& \tilde W(3(k-1))=V(3(k-1))-3(k-1)\\
&=:& z(2j+1),~\mbox{ where } j:=[3(k-1)-1]/2;\\
W^k(1) &=& \tilde W((3(k-1)+1)=V(3(k-1)+1)-[3(k-1)+1]=h\\
W^k(2) &=& \tilde W((3(k-1)+2)=V(3(k-1)+2)-[3(k-1)+2]\\
&=:& z(2j+1),~\mbox{ where } j:=[3(k-1)+1]/2.
\end{eqnarray*}
\end{itemize}
The maximum in (\ref{ape27}) for $k\ge 1$ is provided only by the decisions, with some abuse of notations, represented as  
\begin{eqnarray} \label{ape30}
d_o&:=&(a_2,a_1,a_2)~\mbox{ if $k$ is odd, and}\\
d_e&:=&(a_1,a_2,a_1)~\mbox{  if $k$ is even.} \nonumber
\end{eqnarray}
According to Lemma \ref{apl21}(iii),
$$z(2j+1)+2j+1=(1-\varepsilon)[z(2j-1)+2j-1]+(1+\varepsilon)h+2\varepsilon(j-1),~~j\ge 1,$$
i.e.,
$$z(2j+1)=(1-\varepsilon)z(2(j-1)+1)-2-\varepsilon+(1+\varepsilon)h,$$
and $z(1)=V(1)-1=h-1$. It is obvious that the sequence $\{z(2j+1)\}_{j=0}^\infty$ increases, $z(2j+1)<\left[h+\frac{h-2}{\varepsilon}-1\right]<h$ for all $j\ge 0$ (the rigorous proof is by induction), and $\lim_{j\to\infty} z(2j+1)=h+\frac{h-2}{\varepsilon}-1$. Therefore, the sequence of vectors $\{W^k\}_{k=0}^\infty$ in ${\mathbb R}^M$ is uniformly bounded and
$\lim_{k\to\infty}sp(W^{k+1}-W^k) = 2\left[1-\frac{h-2}{\varepsilon}\right]>0$ as  expected because, if $\lim_{k\to\infty} sp(W^{k+1}-W^k)=0$ then, by the Turnpike Theorem \ref{apt2}, the maximum in (\ref{ape27}) would have been provided by the average optimal decision $d^*$ for large enough $k$.

Condition \ref{apcon2}  is violated because ${\cal Q}(d^*)=\left(
\begin{array}{ccc} 0 & 1 & 0 \\ 0 & 0 & 1 \\ 0 & 1 & 0 \end{array}\right)$, and the state $0$ is here transient.

Condition \ref{apcon3} is violated because the matrix ${\cal Q}(d^*)$ is not aperiodic and, for each $J\ge 1$, for $d^1=d^2=\{d^*,d^*,\ldots,d^*\}$, we have 
$$\eta(d^1,d^2)=\sum_{l=0}^2 \min\{{\cal Q}^{(J,1)}_{1,l},{\cal Q}^{(J,2)}_{2,l}\}=0~\Longrightarrow~\gamma^J=1.$$

Corollary \ref{apcor3} is not applicable because $p_3(a_1)=0$.

The both decision $d_o$ and $d_e$ (\ref{ape30}) are not average optimal (see Lemma \ref{apl22}). At the same time, the switching control strategies $\pi_1:=(d_e,d_0,d_e,\ldots)$ and $\pi_2:=(d_o,d_e,d_o,\ldots)$ are average optimal in the sense that
\begin{equation}\label{ape31}
\lim_{T\to\infty} \frac{1}{T} {\mathbb E}^{\pi_{1,2}}_s\left[\sum_{t=1}^T{\cal R}_{S(t-1)}(d_t)\right]=0,~~~~~s\in{\bf S}.
\end{equation}
(Cf (\ref{ape11}) and see Lemma \ref{apl23} below.) Here ${\mathbb E}^{\pi_{1,2}}_s$ is the mathematical expectation with respect to the strategical measure ${\mathbb P}^{\pi_{1,2}}_s$ on the trajectories $\omega=(s_0=s,s_1,s_2,\ldots)$ in $\bf S$ , which is defined in the standard for MDP way: $S(t,\omega)=s_t$ is the $t$-th component of $\omega$ (the argument $\omega$ is usually omitted), the process $S(t),~t=0,1,\ldots$ is Markov, and ${\mathbb P}^{\pi_{1,2}}_s(S(0)=s)=1$;
$${\mathbb P}^{\pi_{1,2}}_s(S(t)=j|S(t-1)=i)=\left\{
\begin{array}{ll}
{\cal Q}_{i,j}(d_e) & \mbox{ for $\pi_1$, if $t$ is odd}\\
& \mbox{ and for $\pi_2$, if $t$ is even};\\
{\cal Q}_{i,j}(d_o) & \mbox{ for $\pi_1$, if $t$ is even}\\
& \mbox{ and for $\pi_2$, if $t$ is odd}.
\end{array}\right.$$
Similar rigorous constructions are presented in many monographs on MDP: see, e.g., \cite{apb2,apb1}.

\begin{lemma}\label{apl23}
Both for $\pi_1$ and $\pi_2$, the equality (\ref{ape31}) is valid.
\end{lemma}

The situation when a combination of two (ore more) non-optimal control strategies  results in an optimal one, is called `Parrondo's Paradox', see \cite[\S4.2.30]{apb7}.

\section{Conclusion} \label{apis9}

We developed the Turnpike Theory for the undiscounted MDPs and showed that, without appropriated conditions, the turnpike property may be not valid. The presented theory is illustrated by an example of random walk, which seems to be important by its own.

\section{Appendix}\label{apsec4}

\underline{Proof of Lemma \ref{apl13}.} All the elements of the matrix ${\cal Q}(d)$ are obviously non-negative. 

Clearly,
$$\sum_{j=0}^{M-1} {\cal Q}_{0,j}( d)=\sum_{j=0}^{M-1} P_{0,j}( d(0))=1.$$
Suppose $\sum_{j=0}^{M-1} {\cal Q}_{i,j}( d)=1$ for all $i\le l$ for some $l\in\{0,1,\ldots,M-2\}$ and consider $l+1$:
\begin{eqnarray*}
\sum_{j=0}^{M-1} {\cal Q}_{l+1,j}( d)&=& \sum_{j=0}^{l}\sum_{i=0}^l \left[P_{l+1,i}( d(l+1)) {\cal Q}_{i,j}( d)\right]\\
&&+\sum_{j=l+1}^{M-1}\left(\sum_{i=0}^l \left[P_{l+1,i}( d(l+1)) {\cal Q}_{i,j}( d)\right]+P_{l+1,j}( d(l+1))\right)\\
&=& \sum_{i=0}^l\left(\sum_{j=0}^{M-1} {\cal Q}_{i,j}( d)\right) P_{l+1,i}( d(l+1))+\sum_{j=l+1}^{M-1}P_{l+1,j}( d(l+1))\\
&=&\sum_{i=0}^{M-1} P_{l+1,i}( d(l+1))=1.
\end{eqnarray*}
The last  equality is by the induction supposition. 
\hfill$\Box$\vspace{3mm}

\underline{Proof of Proposition \ref{approp1}.} All the statements are well known for the standard model \cite{apb8,apb1}. In the general case, the proofs are word by word identical. We sketch them below.

(a-b) The both operator $U^\beta$ with $\beta\in(0,1)$ and $V\to {\cal R}(d)+\beta{\cal Q}(d)V$ are contractions in the space ${\mathbb R}^M$ with the uniform norm. The maximum in (\ref{ape7}) is provided by some $d^*\in{\bf D}$ according to Assumption \ref{apas2}.

The solutions $V^{\beta *}$ and $V^{\beta d}$ to the equations (\ref{ape7}) and (\ref{ape8}) can be built by iterations like (\ref{ape5}) with $W^0=0$, leading to the sequences $W^{*k}$ and $W^k$ correspondingly.  On each step $W^{* k}\ge W^k$; hence
$$V^{\beta *}=\lim_{k\to\infty} W^{* k}\ge \lim_{k\to\infty} W^k=V^{\beta d}.$$

(c) See \cite[Thm.5.3]{apb8} and also \cite{apb9}.

(d) See \cite[Thm.A.7]{apb1}.

(e) For the existence of a Blackwell optimal decision, see \cite[Thm.4.1]{apb8} which proof remains unchanged in the general case. The same concerns the average optimality of each Blackwell optimal decision: see Thm.5.8(2) and Cor.5.3 of \cite{apb8}. See also \cite[Thm.8.4.5]{apb1}, where the similar statements were proved under the (unneeded) condition that MDP is unichain.
\hfill$\Box$\vspace{3mm}

\underline{Proof of Lemma \ref{apl1}.}  The reasoning is similar to the proof of Theorm 8.5.5 of \cite{apb1}. Below, we provide the key statements.

According to Proposition \ref{approp1}(c), ${\cal Q}^*(\hat d){\cal Q}(\hat d)={\cal Q}^*(\hat d)$. Thus
\begin{eqnarray*}
&&{\cal Q}^*(\hat d){\cal R}(\hat d) = {\cal Q}^*(\hat d)\left[ {\cal R}(\hat d)+{\cal Q}(\hat d)W-W\right]={\cal Q}^*(\hat d)\left[ U\circ W-W\right]\\
&\ge & \min_{s\in{\bf S}}\left\{U\circ W(s)-W(s)\right\}e
\end{eqnarray*}
because ${\cal Q}^*(\hat d)$ is a stochastic matrix.

For any average optimal decision $d^*$, which exists due to Proposition \ref{approp1}(e), we similarly have
\begin{eqnarray*}
&&{\cal Q}^*(\hat d){\cal R}(\hat d)\le{\cal Q}^*(d^*){\cal R}(d^*) = {\cal Q}^*(d^*)\left[ {\cal R}(d^*)+{\cal Q}(d^*)W-W\right]\\
&\le &{\cal Q}^*(d^*)\left[ U\circ W-W\right]
\le  \max_{s\in{\bf S}}\left\{U\circ W(s)-W(s)\right\}e.
\end{eqnarray*}
\hfill$\Box$\vspace{3mm}

\underline{Proof of Lemma \ref{apt1}.} For the standard model Item (a)  follows from Theorem 8.5.2 of \cite{apb1}. For  the general  model $\langle\bf D,{\cal R},{\cal Q}\rangle$, one can repeat the proof of Theorem 8.5.2 \cite{apb1} word by word.

(b) Take $N_1$ such that, for all $n\ge N_1$,
$$sp(U^{nJ}\circ W^1-U^{nJ}\circ W^0)=sp(W^{nJ+1}-W^{nJ})\le\varepsilon;$$
take $N_2\ge N_1$ such that, for all $n\ge N_2$,
$$sp(U^{nJ}\circ W^2-U^{nJ}\circ W^1)=sp(W^{nJ+2}-W^{nJ+1})\le\varepsilon;$$
and so on; \\
take $N_J\ge N_{J-1}\ge\ldots\ge N_1$ such that, for all $n\ge N_J$,
$$sp(U^{nJ}\circ W^J-U^{nJ}\circ W^{J-1})=sp(W^{nJ+J}-W^{nJ+J-1})\le\varepsilon.$$
Now, for $K=N_J$, $sp(W^{k+1}-W^k)\le\varepsilon$ for all $k\ge K$. \hfill$\Box$\vspace{3mm}

\underline{Proof of Corollary \ref{apco2}.} For any two non-negative $M\times M$ matrices $P^1$ and $P^2$ we write $P^1\preceq P^2$ if, for each pair $(s,l)\in{\bf S}^2$, $P^1_{s,l}>0\Longrightarrow P^2_{s,l}>0$.
Since ${\cal Q}_{s,s}(d)>0$ for all $s\in{\bf S}$,  we have
\begin{eqnarray}\label{ape12}
&&P^1\preceq P^1{\cal Q}(d) ~~\mbox{ and}\\
\mbox{if}~P^1\preceq P^2,& \mbox{then}& P^1{\cal Q}(d)\preceq P^2{\cal Q}(d)\label{ape13}
\end{eqnarray}
for all $d\in{\bf D}$ and for all non-negative $M\times M$ matrices $P^1$ and $P^2$.

For each $d\in{\bf D}$, there exists $N_d$ such that ${\cal Q}^n(d)>0$ for all $n\ge N_d$ \cite[Thm.4.1.2]{apb9}. Hence, since the set $\bf D$ is finite, for $N:=\max_{d\in{\bf D}} N_d<\infty$, for all $d\in{\bf D}$, ${\cal Q}^n(d)>0$ for all $n\ge N$. Let $J:=(N-1)|{\bf D}|+1$, where $|{\bf D}|$ is the total number of decisions in $\bf D$.

Fix an arbitrary sequence $\{d_1,d_2,\dots,d_J\}$ of decisions from $\bf D$. Then at least one decision $\hat d$ appears at least $N$ times in this list, say,
$$d_{j_1}=\hat d,~d_{j_2}=\hat d,\ldots, d_{j_N}=\hat d~\mbox{ with } 0<j_1<j_2<\ldots<j_N\le J,$$
and ${\cal Q}^N(\hat d)>0$.

According to (\ref{ape12}),(\ref{ape13}), 
$$I\preceq {\cal Q}(d_1)\preceq {\cal Q}(d_1){\cal Q}(d_2)\preceq\ldots\preceq  {\cal Q}(d_1){\cal Q}(d_2)\ldots{\cal Q}(d_J)={\cal Q}^{(J)},$$
where $I$ is the identical $M\times M$ matrix.

According to (\ref{ape13}), 
$${\cal Q}(\hat d)={\cal Q}(d_{j_1})\preceq{\cal Q}(d_1){\cal Q}(d_2)\ldots P(d_{j_1})$$
because $I\preceq {\cal Q}(d_1){\cal Q}(d_2)\ldots P(d_{j_1-1})$. Similarly
$${\cal Q}^2(\hat d)={\cal Q}(\hat d){\cal Q}(\hat d)\preceq{\cal Q}(d_1){\cal Q}(d_2)\ldots P(d_{j_2})$$
because ${\cal Q}(\hat d)\preceq {\cal Q}(d_1){\cal Q}(d_2)\ldots P(d_{j_1})\preceq{\cal Q}(d_1){\cal Q}(d_2)\ldots P(d_{j_2-1})$.\\
And so on:
$${\cal Q}^N(\hat d)={\cal Q}^{N-1}(\hat d){\cal Q}(\hat d)\preceq{\cal Q}(d_1){\cal Q}(d_2)\ldots P(d_{j_N})$$
because ${\cal Q}^{N-1}(\hat d)\preceq {\cal Q}(d_1){\cal Q}(d_2)\ldots P(d_{j_{N-1}})\preceq{\cal Q}(d_1){\cal Q}(d_2)\ldots P(d_{j_N-1})$.

Therefore, ${\cal Q}^N(\hat d)\preceq {\cal Q}^{(J)}$ meaning that ${\cal Q}^{(J)}>0$, and $\gamma^J<1$ by Lemma \ref{apt3}.\\$~$ \hfill$\Box$

\underline{Proof of Lemma \ref{apt4}.}  All the reasoning is similar to the proof of Corollary 9.4.6 of \cite{apb1}. The milestones are provided below.

(i) For each decision $d\in{\bf D}$, let
$$g(d):={\cal Q}^*(d){\cal R}(d)~\mbox{ and } h(d):=D(d){\cal R}(d)$$
be the so called gain and bias of the decision $d$. Then
$$g(d)={\cal Q}(d)g(d)~\mbox{ and } g(d)+(I-{\cal Q}(d))h(d)={\cal R}(d).$$
(See \cite[Thm.8.2.6]{apb1}.)

(ii) There exists a solution $(g,Y)\in {\mathbb R}^M\times{\mathbb R}^M$ to the (so called `modified') optimality equation
\begin{equation}\label{ape14}
g=\max_{d\in{\bf D}}\{{\cal Q}(d)g\};~~~~~ g+Y=\max_{d\in{\bf D}}\{{\cal R}(d)+{\cal Q}(d)Y\},
\end{equation}
for which necessarily $g=g^*:=\max_{d\in{\bf D}}{\cal Q}^*(d){\cal R}(d)$ is the maximal gain \cite[Cor.5.4]{apb8}, \cite[Thm.9.1.2, Cor.9.1.5]{apb1}. If a decision $d^*$ provides the both maxima, then $d^*$ is average optimal \cite[Thm.9.1.7]{apb1}.

(iii) $\displaystyle \lim_{k\to\infty}\frac{1}{k} W^k=g^*$.\\
See \cite[Thm.9.4.1]{apb1} and also \cite[Lemma 5.5]{apb8}.

(iv) Suppose, for every average optimal decision $d^*\in{\bf D}^*$, the matrix ${\cal Q}(d^*)$ is aperiodic. Then there exists the limit
$$\lim_{k\to\infty}[W^k-kg^*]$$
 \cite[Thm.9.4.4]{apb1}.
 
(v) Therefore,
$${\bf 0}=\lim_{k\to\infty}\left[ [W^{k+1}-(k+1)g^*]-[W^k-k g^*]\right]=\lim_{k\to\infty}\left[(W^{k+1}-W^k)-g^*\right],$$
where ${\bf 0}$ is the zero vector in ${\mathbb R}^M$.

(vi) Under the assumptions of Lemma \ref{apt4},
$$\lim_{k\to\infty} (W^{k+1}-W^k)=G^*e\Longrightarrow\lim_{k\to\infty}~sp(W^{k+1}-W^k)=0.$$

All the statements (i)-(vi) formulated above were proved in \cite{apb8,apb1} for the standard model; all the steps of the proofs remain valid for the general case.\linebreak
$~$\hfill$\Box$

\underline{Proof of Lemma \ref{apl8}.} The case of $i=-M,-M+1,\ldots, -1$ is obvious.

For $i\ge 0$,
\begin{eqnarray*}
\tilde W(i) &=& \max_{a\in{\bf A}}\left\{ R^a+\sum_{m=1}^M[\tilde W(i-m)+c_*(i-m)] p_m(a)\right\}-c_*i\\
&=& \max_{a\in{\bf A}} \left\{ R^a-c_* L^a+\sum_{m=1}^M \tilde W(i-m)p_m(a)\right\}.
\end{eqnarray*}
Equalities (\ref{ape23}) are proved, and the maxima in (\ref{ape22}) and in (\ref{ape23}) are provided by the same values of $a\in{\bf A}$.

Finally, keeping in mind that 
\begin{itemize}
\item[$\bullet$] $|\tilde W(i)|\le |c_*|M$ for $i<0$,
\item[$\bullet$]  $\frac{R^a}{L^a}-c_*\le 0$ for all $a\in{\bf A}$, and
\item[$\bullet$]  $\frac{R^a}{L^a}-c_*= 0$ for  $a\in{\bf A}_*\ne\emptyset$,
\end{itemize}
it is easy to prove by induction that $|\tilde W(i)|\le |c_*|M$ for all $i=0,1,2,\ldots$.
  \hfill $\Box$\vspace{3mm}
  
\underline{Proof of Lemma \ref{apl11}.}  Suppose $d^*\in{\bf D}^*$, that is, 
$${\cal Q}^*(d^*){\cal R}(d^*)=\max_{d\in{\bf D}}{\cal Q}^*(d){\cal R}(d)={\bf 0}.$$ Matrix ${\cal Q}^*(d^*)$ in the block representation (`standard form') has only the blocks corresponding to the closed communicating classes, and those blocks are strictly positive \cite[\S5.3.1]{apb8}. Thus, if ${\cal R}(d^*)\ne{\bf 0}$ then the non-positive vector ${\cal Q}^*(d^*){\cal R}(d^*)$ contains negative elements, and $d^*\notin{\bf D}^*$. The obtained contradiction completes the proof.   \hfill $\Box$\vspace{3mm}

\underline{Proof of Corollary \ref{apcor3}.} Let us show that, for each $d\in{\bf D}$, for any two states $i,j\in{\bf S}$, ${\cal Q}_{i,j}(d)>0$ provided $P_{i,j}(d(i))>0$.

If $j\ge i$ then this statement follows directly from the definition of the matrix ${\cal Q}(d)$:
$${\cal Q}_{i,j}(d)\ge P_{i,j}(d(i))>0.$$

Suppose $j<i$. Then, again using the definition of he matrix ${\cal Q}(d)$, we have
$${\cal Q}_{i,j}(d)\ge P_{i,j}(d(i)){\cal Q}_{j,j}(d).$$
Since ${\cal Q}_{j,j}(d)\ge P_{j,j}(d(j))=p_M(d(j))>0$, we obtain the required inequality ${\cal Q}_{i,j}(d)>0$, if $P_{i,j}(d(i))>0$.

Now, for any two states $i,j\in{\bf S}$, for the path $i_0=i\to i_1\to\ldots\to i_N=j$ in $\bf S$, we have
$${\cal Q}_{i_0,i_1}(d){\cal Q}_{i_1,i_2}(d)\ldots{\cal Q}_{i_{N-1},i_N}(d)>0$$
for each $d\in{\bf D}$. To put it different, for each $d\in{\bf D}$, the matrix ${\cal Q}(d)$ corresponds to the irreducible Makov chain which is certainly aperiodic because ${\cal Q}_{i,i}(d)\ge P_{i,i}(d(i))>0$ for all $i\in{\bf S}$.

Since Conditions \ref{apcon2} and \ref{apcon3} are satisfied, the statement of Corollary follows from Theorem \ref{apt6}.
 \hfill $\Box$\vspace{3mm}

\underline{Proof of Lemma \ref{apl21}.} 
When $j=1$, Items (i) and (iii) are valid by the preliminary calculations, and Item (ii) comes from the following:
$$2\varepsilon(j-1)+(1+\varepsilon)h-2-\varepsilon V(2j-1)=(1+\varepsilon)h-2-\varepsilon h=h-2>0.$$

Suppose statements (i), (ii) and (iii) hold for some $j\ge 1$ and consider $j+1$.

(i) For $i=2(j+1)$, using the induction supposition, we estimate the difference
\begin{eqnarray*}
&& 2+V(2j)-[h+(1-\varepsilon)V(2j)+\varepsilon V(2j-1)]\\
&>& 2+2j+h-h-(1-\varepsilon)(2j+h)-\varepsilon\frac{2\varepsilon(j-1)+(1+\varepsilon)h-2}{\varepsilon}\\
&=& 4+2\varepsilon-2h=2[2+\varepsilon-h]>0.
\end{eqnarray*}
The  inequality is according to statement (ii) at $j$.

Thus, $V(2(j+1))=2(j+1)+h$, and the maximum in (\ref{ape21}) is provided only by $a_1$.

(ii)
\begin{eqnarray*}
&& V(2j+1)=(1-\varepsilon)V(2j-1)+(1+\varepsilon)h+2\varepsilon(j-1)\\
&<& (1-\varepsilon)\frac{2\varepsilon(j-1)+(1+\varepsilon)h-2}{\varepsilon}+(1+\varepsilon)h+2\varepsilon(j-1)\\
&=&\frac{2\varepsilon j+(1+\varepsilon)h-2}{\varepsilon},
\end{eqnarray*}
so that statement (ii) is valid for $j+1$.

(iii) For $i=2(j+1)+1$, using the induction supposition, we estimate the difference
\begin{eqnarray*}
&& h+(1-\varepsilon)V(2j+1)+\varepsilon V(2j) -[2+V(2j+1)]\\
&=& h-\varepsilon V(2j+1)+\varepsilon[2j+h]-2\\
&>& h(1+\varepsilon)-[2\varepsilon j+h-2+\varepsilon h]+2\varepsilon j-2=0,
\end{eqnarray*}
where the inequality is by the proved above Item (ii) for $j+1$. Recall also that  $V(2j)=2j+h$. Therefore,
$$V(2(j+1)+1)=h+(1-\varepsilon) V(2j+1)+\varepsilon[2j+h),$$
and we see that statement (iii) is  valid for $j+1$ and the maximum in (\ref{ape21}) is provided only by $a_2$.
 \hfill $\Box$\vspace{3mm}
 
\underline{Proof of Lemma \ref{apl22}.} One should consider all the different decisions $d:{\bf S}\to{\bf A}$. Since the reasoning is similar in all the cases, below we investigate only the decision $d(0)=a_1$, $d(1)=a_1$, $d(2)=a_2$. Here
 $${\cal R}(d)=(0,0, h-(2+\varepsilon))^T;~~~~~
{\cal Q}(d) =
\left(\begin{array}{lll}
0 & 1 & 0 \\ 0 & 0 & 1 \\ 0 & 1-\varepsilon & \varepsilon
\end{array}\right).$$
Recall that the model $\langle{\bf D},{\cal R},{\cal Q}\rangle$ is as in the version (ii) described in Section \ref{apsec2}; $r_s(a):=L^a\left(\frac{R^a}{L^a}-c_*\right)$ and the matrix $P(a)$ is given by the equation (\ref{ape26}). Remember also that $h-(2+\varepsilon)<0$.
Clearly, the states $1$ and $2$ form the closed communicating class, the state $0$ is transient, and the components ${\cal Q}^*_{0,2}(d),{\cal Q}^*_{1,2}(d),{\cal Q}^*_{2,2}(d)$ are strictly positive, so that
${\cal Q}^*(d){\cal R}(d)<0$,
and the decision $d$ is not average optimal.

For all other non-trivial decisions $d$, one can similarly show that the vector ${\cal Q}^*(d){\cal R}(d)$ contains strictly negative components. \hfill $\Box$\vspace{3mm}

\underline{Proof of Lemma \ref{apl23}.} The reasoning is similar for $\pi_1$ and $\pi_2$, so we consider only $\pi_1=(d_e,d_0,d_e,\ldots)$. Obviously, it is sufficient to show that the sequence of vectors in ${\mathbb R}^M$
$$R^T:=\left\{{\mathbb E}^{\pi_1}_s \left[\sum_{t=3}^T{\cal R}_{S(t-1)}(d_t)\right],~s\in{\bf S}\right\},~~~T\ge 3,$$
is uniformly bounded.

Suppose an odd value of $T\ge 3$ is fixed and consider the sequence of vectors $\{\tilde R^l_i,~i\in{\bf S}\}_{l=2}^T$ in ${\mathbb R}^M$ defined as follows:
$$\tilde R^l_i:={\mathbb E}^{\pi_1}_s\left[\sum_{t=T-l+3}^T{\cal R}_{S(t-1)}(d_t)+W^2(S(T))|S(T-l+2)=i\right].$$
These vectors are $s$-independent due to the Markov property of the process $\{S(t),~t=0,1,\ldots\}$. Since the maximum in (\ref{ape27}) is provided by $d_e$ ($d_o$) if $k$ is even (odd),
\begin{eqnarray*}
\tilde R^2&=&W^2;\\
\tilde R^3&=&{\cal R}(d_e)+{\cal Q}(d_e)\tilde R^2=W^3;\\
\tilde R^4&=&{\cal R}(d_o)+{\cal Q}(d_o)\tilde R^3=W^4;~~\ldots,\\
\tilde R^T&=&\left\{{\mathbb E}^{\pi_1}_s\left[\sum_{t=3}^T{\cal R}_{S(t-1)}(d_t)+W^2(S(T))|S(2)=i\right],~i\in{\bf S}\right\}\\
&=&{\cal R}(d_e)+{\cal Q}(d_e)\tilde R^{T-1}=W^T.
\end{eqnarray*}
The sequence $\{W^T\}_{T=0}^\infty$ is uniformly bounded by Lemma \ref{apl8}; so the sequence of vectors
\begin{eqnarray*}
R^T &=&\left\{ {\mathbb E}^{\pi_1}_s \left[{\mathbb E}^{\pi_1}_s\left[\sum_{t=3}^T{\cal R}_{S(t-1)}(d_t)+W^2(S(T))|S(2)\right]\right]\right.\\
&&\left.\vphantom{\sum_{t=3^T}}-{\mathbb E}^{\pi_1}_s \left[{\mathbb E}^{\pi_1}_s\left[W^2(S(T))|S(2)\right]\right],~s\in{\bf S}\right\}\\
&=&\left\{ {\mathbb E}^{\pi_1}_s\left[W^T(S(2))\right]-{\mathbb E}^{\pi_1}_s\left[W^2(S(T))\right],~s\in{\bf S}\right\},~~~T=3,5,7,\ldots
\end{eqnarray*}
is uniformly bounded.

For even values of $T$, the sequence 
$$R^T:=\left\{{\mathbb E}^{\pi_1}_s \left[\sum_{t=3}^{T-1}{\cal R}_{S(t-1)}(d_t)\right]+{\mathbb E}^{\pi_1}_s\left[{\cal R}_{S(T-1)}(d_T)\right],~s\in{\bf S}\right\},~~~T=4,6,8,\ldots$$
is also uniformly bounded, and the proof is completed.
\hfill $\Box$\vspace{3mm}


\end{document}